\documentclass{article}
\usepackage{amssymb}
\usepackage{amsfonts}
\usepackage{amsmath}

\numberwithin{theorem}{section}
\numberwithin{equation}{section}

\begin{document}

\title{About intrinsic Finsler connections
for the homogeneous lift to the Osculator Bundle of a Finsler metric}
\author{Alexandru Oana\footnotemark[1]  \;}
\date{}
\maketitle

\begin{abstract}
In this article we present a study of the subspaces of the manifold OscM, the
total space of the osculator bundle of a real manifold M. We obtain the
induced connections of the canonical metrical N-linear connection
determined by the homogeneous prolongation of a Finsler metric to the
manifold OscM. We present the relation between the induced and intrinsic
geometric objects of the associated osculator submanifold.
\end{abstract}

\footnotetext[1]{%
Department of Mathematics and Informatics, University Transilvania of Bra%
\c{s}ov, 50 Iuliu Maniu Blvd., 500091 Bra\c{s}ov, Romania, E-mail: 
alexandru.oana@unitbv.ro} \setcounter{footnote}{3}



\begin{center}
\bigskip {\LARGE Introduction}
\end{center}

The Sasaki $N$-prolongation $\mathbb{G}$ to the osculator bundle without
the null section $\widetilde{OscM}=OscM\backslash \left\{ 0\right\} $ of a
Finslerian metric $g_{ab}$ on the manifold $M$ given by%
\begin{equation*}
\mathbb{G}=g_{ab}\left( x,y\right) dx^{a}\otimes
dx^{b}+g_{ab}\left( x,y\right)
\delta y^a\otimes \delta y^b
\end{equation*}%
is a Riemannian structure on $\widetilde{OscM},$ which depends only on
the metric $g_{ab}.$ 

The tensor $\mathbb{G}$ is not invariant with respect to the homothetis on
the fibres of $\widetilde{OscM}$, because $\mathbb{G}$ is not homogeneous with respect
to the variable $y^{a}.$

In this paper, we use a new kind of prolongation $\mathbb{\mathring{G}}$ to $\widetilde{OscM}$, \cite{at3}, which depends
only on the metric $g_{ab}.$ Thus, $\mathbb{\mathring{G}}$ determines on the
manifold $\widetilde{OscM}$ a Riemannian structure which is 0-homogeneous on
the fibres of $OscM.$

Some geometrical properties of $\mathbb{\mathring{G}}$ are studied: the
canonical metrical $N$-linear connection, the induced linear connections etc.

\section{Preliminaries}

Let us consider $F^{n}=\left( M,F\right) $ a Finsler space $\left( \cite%
{mian2}\right) $, and $F:TM=OscM\rightarrow \mathbb{R}$ the fundamental
function. $F$ is a $C^{\infty }$ function on the manifold $OscM$ and it is continuous on the null section
of the projection $\pi :OscM\rightarrow M.$ The fundamental tensor on $F^{n}$
is

\begin{equation*}
g_{ab}\left( x,y\right) =\dfrac{1}{2}\dfrac{\partial ^{2}F^{2}}{\partial
y^{a}\partial y^{b}},\;\;\forall \left( x,y\right) \in OscM.
\end{equation*}

The lagrangian $F^{2}\left( x,y\right) $ determines the canonical spray $%
S=y^{a}\dfrac{\partial }{\partial x^{a}}-2G^{a}\dfrac{\partial }{\partial
y^{a}}$ with the coefficients $G^{a}=\dfrac{1}{2}\gamma \overset{}{_{bc}^{a}}%
\left( x,y\right) y^{b}y^{c},$ where $\gamma \overset{}{_{bc}^{a}}\left(
x,y\right) $ are the Christoffels symbols of the metric tensor $g_{ab}\left(
x,y\right) .$ The Cartan nonlinear connection $N$ of the space $F^{n}$%
\textbf{\ }has the coefficients

\begin{equation}
N\overset{}{^{a}}_{b}=\dfrac{\partial G^{a}}{\partial y^{b}}.  \label{cncvar}
\end{equation}%
$N$ determines a distribution on the manifold $\widetilde{OscM},$ (\cite%
{mian2},\cite{mian1}), which is supplementary to the vertical distribution $%
V.$ We have the next decomposition

\begin{equation}
T_{w}\widetilde{OscM}=N_{w}\oplus V_{w},\forall w=\left( x,y\right) \in 
\widetilde{Osc}M.  \label{desc}
\end{equation}

The adapted basis of this decomposition is $\left\{ \dfrac{\delta }{\delta
x^{a}},\dfrac{\partial }{\partial y^{a}}\right\} $ , $\left( a=1,..,n\right) 
$ and its dual basis is $\left( dx^{a},\delta y^{a}\right) ,$ where

\begin{equation}
\left\{ 
\begin{array}{l}
\dfrac{\delta }{\delta x^{a}}=\dfrac{\partial }{\partial x^{a}}-N\overset{}{%
^{b}}_{a}\dfrac{\delta }{\delta y^{b}}, \\ 
\\ 
\dfrac{\partial }{\partial y^{a}}=\qquad \quad \quad \qquad \dfrac{\partial 
}{\partial y^{a}}%
\end{array}%
\right.  \label{bz_ad}
\end{equation}

and

\begin{equation}
\left\{ 
\begin{array}{l}
dx^{a}=%
\begin{array}{lll}
&  & 
\end{array}%
dx^{a}, \\ 
\delta y^{a}=dy^{a}+N\overset{}{^{a}}_{b}dx^{b}.%
\end{array}%
\right.  \tag{1.1.7}  \label{bz_du}
\end{equation}

We use the next notations: 
\begin{equation*}
\delta _{a}=\dfrac{\delta }{\delta x^{a}},\;\;\dot{\partial}_{1a}=\dfrac{%
\partial }{\partial y^{a}}.
\end{equation*}

The fundamental tensor $g_{ab}$ determines on the manifold $\widetilde{OscM}$
the homogeneous N-lift $\overset{0}{\mathbb{G}},$\cite{Mir+Shi+Sab}$,$

\begin{equation}
\overset{0}{\mathbb{G}}=g_{ab}\left( x,y\right) dx^{a}\otimes
dx^{b}+h_{ab}\left( x,y\right) \delta y^{a}\otimes \delta y^{b},
\label{lift_om}
\end{equation}%
where

\begin{equation}
h_{ab}\left( x,y\right) =\frac{p^{2}}{\left\Vert y\right\Vert ^{2}}%
g_{ab}\left( x,y\right) ,  \tag{1.1.9}  \label{comp_lift}
\end{equation}

\begin{equation*}
\left\Vert y\right\Vert ^{2}=g_{ab}\left( x,y\right) y^{a}y^{b}.
\end{equation*}

This is homogeneous with respect to $y,$ and $p$ is a positive constant
required by applications in order that the physical dimensions of the terms
of $\mathbb{\mathring{G}}$ be the same.

Let $\check{M}$ be a real, m-dimensional manifold, immersed in $M$ through
the immersion $i:\check{M}\rightarrow M$. Localy, $i$ can be given in the
form

\begin{equation*}
\begin{array}{lcl}
x^{a}=x^{a}\left( u^{1},...,u^{m}\right) , &  & rank\left\Vert \dfrac{%
\partial x^{a}}{\partial u^{\alpha }}\right\Vert =m.%
\end{array}%
\end{equation*}

The indices $a,b,c,$....run over the set $\left\{ 1,...,n\right\} $ and $%
\alpha ,\beta ,\gamma ,...$ run on the set $\left\{ 1,...,m\right\} .$ We
assume $1<m<n$. We take the immersed submanifold $Osc\check{M}$ $\ $of the
manifold $OscM,$ by the immersion Osc$i:Osc\check{M}\rightarrow OscM.$ The
parametric equations of the submanifold $Osc\check{M}$ are

\begin{equation}
\left\{ 
\begin{array}{l}
x^{a}=x^{a}\left( u^{1},...,u^{m}\right) ,rang\left\Vert \dfrac{\partial
x^{a}}{\partial u^{\alpha }}\right\Vert =m \\ 
\\ 
y^{a}=\dfrac{\partial x^{a}}{\partial u^{\alpha }}v^{\alpha }.%
\end{array}%
\right.  \label{ec_subvar}
\end{equation}

The restriction of the fundamental function $F$ to the submanifold $%
\widetilde{Osc\check{M}}$ is

\begin{equation*}
\check{F}\left( u,v\right) =F\left( x\left( u\right) ,y\left( u,v\right)
\right)
\end{equation*}%
and we call $\check{F}^{m}=\left( \check{M},\check{F}\right) $ the \textbf{%
induced Finsler subspaces} of $F^{n}$ and $\check{F}$ \ the \textbf{induced
fundamental function}.

Let $B_{\alpha }^{a}(u)=\dfrac{\partial x^{a}}{\partial u^{\alpha }}$ and $%
g_{\alpha \beta }$ the induced fundamental tensor, 
\begin{equation}
g_{\alpha \beta }\left( u,v\right) =g_{ab}\left( x\left( u\right) ,y\left(
u,v\right) \right) \overset{}{B_{\alpha }^{a}}\overset{}{B_{\beta }^{b}}.
\label{tens_fund_subvar}
\end{equation}%
We obtain a system of d-vectors $\left\{ B_{\alpha }^{a},B_{\bar{\alpha}%
}^{a}\right\} $ wich determines a moving frame $\mathcal{R=}\left\{ \left(
u,v\right) ;B_{\alpha }^{a}\left( u\right) ,B_{\bar{\alpha}}^{a}\left(
u,v\right) \right\} $ in $OscM$ along to the submanifold $Osc\check{M}.$

Its dual frame will be denoted by $\mathcal{R}^{\ast }\mathcal{=}\left\{
B_{a}^{\alpha }\left( u,v\right) ,B_{a}^{\bar{\alpha}}\left( u,v\right)
\right\} .$ This is also defined on an open set $\check{\pi}^{-1}\left( 
\check{U}\right) \subset Osc\check{M},$ $\check{U}$ being a domain of a
local chart on the submanifold $\check{M}.$

The conditions of duality are given by:%
\begin{equation*}
B_{\beta }^{a}B_{a}^{\alpha }=\delta _{\beta }^{\alpha },\quad B_{\beta
}^{a}B_{a}^{\bar{\alpha}}=0,\quad B_{a}^{\alpha }B_{\bar{\beta}}^{a}=0,\quad
B_{a}^{\bar{\alpha}}B_{\bar{\beta}}^{a}=\delta _{\bar{\beta}}^{\bar{\alpha}}
\end{equation*}%
\begin{equation*}
\begin{array}{l}
B_{\alpha }^{a}B_{b}^{\alpha }+B_{\bar{\alpha}}^{a}B_{b}^{\bar{\alpha}%
}=\delta _{b}^{a}.%
\end{array}%
\end{equation*}

The restriction of the of the nonlinear connection N to $\widetilde{Osc\check{M}}$
uniquely determines an induced nonlinear connection $\check{N}$ on $\widetilde{Osc\check{M}}$

\begin{equation}
\check{N}\overset{}{^{\alpha }}\overset{}{_{\beta }}=B_{a}^{\alpha }\left(
B_{0\beta }^{a}+N\overset{}{^{a}}\overset{}{_{b}}B_{\beta }^{b}\right) .
\label{cnc_subvar}
\end{equation}

The cobasis $\left( dx^{i},\delta y^{a}\right) $\ restricted to $Osc\check{M}
$ is uniquely represented in the moving frame $\mathcal{R}$\ in the
following form:%
\begin{equation}
\left\{ 
\begin{array}{l}
dx^{a}=B_{\beta }^{a}du^{\beta } \\ 
\\ 
\delta y^{a}=B_{\alpha }^{a}\delta v^{\alpha }+B_{\bar{\alpha}}^{a}K\overset{%
}{_{\beta }^{\bar{\alpha}}}du^{\beta }%
\end{array}%
\right.  \label{1.4}
\end{equation}%
where%
\begin{equation*}
K\overset{}{_{\beta }^{\bar{\alpha}}}=B_{a}^{\bar{\alpha}}\left( B_{0\beta
}^{a}+M\overset{}{_{b}^{a}}B_{\beta }^{b}\right) ,\text{ }B_{0\beta
}^{a}=B_{\alpha \beta }^{a}v^{a}.
\end{equation*}

A linear connection $D$ on the manifold $OscM$ is called \textbf{metrical
N-linear connection} with respect to $\mathbb{\mathring{G}},$ if $D\mathbb{%
\mathring{G}=}0$ and $D$ preserves by parallelism the distributions N and V.
The coefficients of the N-linear connections $D\Gamma \left( N\right) $ will
be denoted with $\left( \underset{\left( 00\right) }{\overset{H}{L}}\overset{%
}{_{bc}^{a}},\underset{\left( 10\right) }{\overset{V}{L}}\overset{}{_{bc}^{a}%
},\underset{\left( 01\right) }{\overset{H}{C}}\overset{}{_{bc}^{a}},\underset%
{\left( 11\right) }{\overset{V}{C}}\overset{}{_{bc}^{a}}\right) .$ \newline
\textbf{Theorem 1.1}(\cite{Mir+Shi+Sab}) \textit{There exist metrical} $N$-%
\textit{linear connections} $D\Gamma \left( N\right) $ \textit{on} $%
\widetilde{OscM},$ \textit{with respect to the homogeneous prolongation} $%
\mathbb{\mathring{G}}$, \textit{wich depend only on the metric} $%
g_{ab}\left( x,y\right) .$ \textit{One of these connections has \newline
the "horizontal" coefficients}%
\begin{equation}
\begin{array}{l}
\underset{\left( 00\right) }{\overset{H}{L}}\overset{}{_{bc}^{a}}=\dfrac{1}{2%
}g^{ad}\left( \delta _{b}g_{dc}+\delta _{c}g_{bd}-\delta _{d}g_{bc}\right)
\\ 
\\ 
\underset{\left( 10\right) }{\overset{V}{L}}\overset{}{_{bc}^{a}}=\dfrac{1}{2%
}h^{ad}\left( \delta _{b}h_{dc}+\delta _{c}h_{bd}-\delta _{d}h_{bc}\right)%
\end{array}
\label{h_coef_var}
\end{equation}%
\textit{and the "vertical" coefficients}:%
\begin{equation}
\begin{array}{l}
\underset{\left( 01\right) }{\overset{H}{C}}\overset{}{_{bc}^{a}}=\dfrac{1}{2%
}g^{ad}\left( \dot{\partial}_{b}g_{dc}+\dot{\partial}_{c}g_{bd}-\dot{\partial%
}_{d}g_{bc}\right) \\ 
\\ 
\underset{\left( 11\right) }{\overset{V}{C}}\overset{}{_{bc}^{a}}=\dfrac{1}{2%
}h^{ad}\left( \dot{\partial}_{b}h_{dc}+\dot{\partial}_{c}h_{bd}-\dot{\partial%
}_{d}h_{bc}\right) .%
\end{array}
\label{v_coef_var}
\end{equation}%
It is called the \textbf{Cartan metrical N-linear connection}. This linear
connection will be used throughout this paper.

For this N-linear connection, we have the operators $\overset{H}{D}$ and $%
\overset{V}{D}$ which are given by the following relations%
\begin{equation}
\begin{array}{l}
\overset{H}{D}X^{a}=dX^{a}+\overset{H}{\omega _{b}^{a}}X^{b} \\ 
\\ 
\overset{V}{D}X^{a}=dX^{a}+\overset{V}{\omega _{b}^{a}}X^{b}%
\end{array}%
\forall X\in \mathcal{F}\left( \widetilde{Osc}M\right) .  \label{dif_var}
\end{equation}%
We call these operators the \textbf{horizontal} and \textbf{vertical
covariant differentials}. The 1-forms which define these operators will be
called the \textbf{horizontal} and \textbf{vertical 1-form}, where%
\begin{equation}
\begin{array}{l}
\overset{H}{\omega }\underset{}{_{b}^{a}}=\underset{\left( 00\right) }{%
\overset{H}{L}}\overset{}{_{bc}^{a}}dx^{c}+\underset{\left( 01\right) }{%
\overset{H}{C}}\overset{}{_{bc}^{a}}\delta y^{c} \\ 
\\ 
\overset{V}{\omega _{b}^{a}}=\underset{\left( 10\right) }{\overset{V}{L}}%
\overset{}{_{bc}^{a}}dx^{c}+\underset{\left( 11\right) }{\overset{V}{C}}%
\overset{}{_{bc}^{a}}\delta y^{c}.%
\end{array}
\label{1_form_var}
\end{equation}%
We have \newline
\textbf{Theorem 1.2} \textit{The} \textit{d-tensors of torsion of the Cartan
metrical N-linear connection }$D$\textit{\ have the next expresions:}%
\begin{equation}
\begin{array}{lll}
\underset{\left( 00\right) }{\overset{H}{T}}\overset{}{_{bc}^{a}}=\underset{%
\left( 00\right) }{\overset{H}{L}}\overset{}{_{bc}^{a}}-\underset{\left(
00\right) }{\overset{H}{L}}\overset{}{_{cb}^{a}}, &  & \underset{\left(
01\right) }{\overset{V}{T}}\overset{}{_{bc}^{a}}=R\overset{}{_{bc}^{a}}, \\ 
&  &  \\ 
\underset{\left( 10\right) }{\overset{H}{P}}\overset{}{_{bc}^{a}}=\underset{%
\left( 01\right) }{\overset{H}{C}}\overset{}{_{bc}^{a}}, &  & \underset{%
\left( 11\right) }{\overset{V}{P}}\overset{}{_{bc}^{a}}=\underset{\left(
11\right) }{B}\overset{}{_{bc}^{a}}-\underset{\left( 10\right) }{\overset{V}{%
L}}\overset{}{_{cb}^{a}} \\ 
&  &  \\ 
\underset{\left( 11\right) }{\overset{V}{S}}\overset{}{_{bc}^{a}}=\underset{%
\left( 11\right) }{\overset{V}{C}}\overset{}{_{bc}^{a}}-\underset{\left(
11\right) }{\overset{V}{C}}\overset{}{_{cb}^{a}}. &  & 
\end{array}
\label{tors_var}
\end{equation}%
\textbf{\bigskip Theorem 1.3} \textit{The Cartan metrical N-linear
connection }$D$\textit{\ has, in the adapted bases }$\left\{ \delta _{a},\dot{%
\partial}_{1a}\right\} $\textit{,\ the following d-tensors of curvature \newline
"horizontals"}%
\begin{equation}
\begin{array}{llllllllll}
\underset{\left( 00\right) }{\overset{H}{R}}\underset{}{_{b}}\overset{}{^{a}}%
\underset{}{_{cd}} & = & \delta _{d}\underset{\left( 00\right) }{\overset{H}{%
L}}\overset{}{_{bc}^{a}} & - & \delta _{c}\underset{\left( 00\right) }{%
\overset{H}{L}}\overset{}{_{bd}^{a}} & + & \underset{\left( 00\right) }{%
\overset{H}{L}}\overset{}{_{bc}^{f}}\underset{\left( 00\right) }{\overset{H}{%
L}}\overset{}{_{fd}^{a}} & - & \underset{\left( 00\right) }{\overset{H}{L}}%
\overset{}{_{bd}^{f}}\underset{\left( 00\right) }{\overset{H}{L}}\overset{}{%
_{fc}^{a}} & + \\ 
&  &  &  &  &  &  &  &  &  \\ 
&  &  &  &  & + & \underset{\left( 01\right) }{\overset{H}{C}}\overset{}{%
_{bf}^{a}}R\overset{}{_{cd}^{f}}, &  &  &  \\ 
&  &  &  &  &  &  &  &  &  \\ 
\underset{\left( 10\right) }{\overset{H}{P}}\underset{}{_{b}}\overset{}{^{a}}%
\underset{}{_{cd}} & = & \dot{\partial}_{1d}\underset{\left( 00\right) }{%
\overset{H}{L}}\overset{}{_{bc}^{a}} & - & \underset{\left( 01\right) }{%
\overset{H}{C}}\overset{}{_{bd\mid 0c}^{a}} & + & \underset{\left( 01\right) 
}{\overset{H}{C}}\overset{}{_{bf}^{a}}\underset{\left( 11\right) }{\overset{H%
}{P}}\overset{}{_{cd}^{f}}, &  &  &  \\ 
&  &  &  &  &  &  &  &  &  \\ 
\underset{\left( 10\right) }{\overset{H}{S}}\underset{}{_{b}}\overset{}{^{a}}%
\underset{}{_{cd}} & = & \dot{\partial}_{1d}\underset{\left( 01\right) }{%
\overset{H}{C}}\overset{}{_{bc}^{a}} & - & \dot{\partial}_{1c}\underset{%
\left( 01\right) }{\overset{H}{C}}\overset{}{_{bd}^{a}} & + & \underset{%
\left( 01\right) }{\overset{H}{C}}\overset{}{_{bc}^{f}}\underset{\left(
01\right) }{\overset{H}{C}}\overset{}{_{fd}^{a}} & - & \underset{\left(
01\right) }{\overset{H}{C}}\overset{}{_{bd}^{f}}\underset{\left( 01\right) }{%
\overset{H}{C}}\overset{}{_{fc}^{a}}, & 
\end{array}
\label{curb_var_h}
\end{equation}%
\newline
\textit{and the "verticals"}%
\begin{equation}
\begin{array}{llllllllll}
\underset{\left( 01\right) }{\overset{V}{R}}\underset{}{_{b}}\overset{}{^{a}}%
\underset{}{_{cd}} & = & \delta _{d}\underset{\left( 10\right) }{\overset{V}{%
L}}\overset{}{_{bc}^{a}} & - & \delta _{c}\underset{\left( 10\right) }{%
\overset{V}{L}}\overset{}{_{bd}^{a}} & + & \underset{\left( 10\right) }{%
\overset{V}{L}}\overset{}{_{bc}^{f}}\underset{\left( 10\right) }{\overset{V}{%
L}}\overset{}{_{fd}^{a}} & - & \underset{\left( 10\right) }{\overset{V}{L}}%
\overset{}{_{bd}^{f}}\underset{\left( 10\right) }{\overset{V}{L}}\overset{}{%
_{fc}^{a}} & + \\ 
&  &  &  &  &  &  &  &  &  \\ 
&  &  &  &  & + & \underset{\left( 11\right) }{\overset{V}{C}}\overset{}{%
_{bf}^{a}}R\overset{}{_{cd}^{f}}, &  &  &  \\ 
&  &  &  &  &  &  &  &  &  \\ 
\underset{\left( 11\right) }{\overset{V}{P}}\underset{}{_{b}}\overset{}{^{a}}%
\underset{}{_{cd}} & = & \dot{\partial}_{1d}\underset{\left( 10\right) }{%
\overset{V}{L}}\overset{}{_{bc}^{a}} & - & \underset{\left( 11\right) }{%
\overset{V}{C}}\overset{}{_{bd\mid 1c}^{a}} & + & \underset{\left( 11\right) 
}{\overset{V}{C}}\overset{}{_{bf}^{a}}\underset{\left( 11\right) }{\overset{V%
}{P}}\overset{}{_{cd}^{f}}, &  &  &  \\ 
&  &  &  &  &  &  &  &  &  \\ 
\underset{\left( 11\right) }{\overset{V}{S}}\underset{}{_{b}}\overset{}{^{a}}%
\underset{}{_{cd}} & = & \dot{\partial}_{1d}\underset{\left( 11\right) }{%
\overset{V}{C}}\overset{}{_{bc}^{a}} & - & \dot{\partial}_{1c}\underset{%
\left( 11\right) }{\overset{V}{C}}\overset{}{_{bd}^{a}} & + & \underset{%
\left( 11\right) }{\overset{V}{C}}\overset{}{_{bc}^{f}}\underset{\left(
11\right) }{\overset{V}{C}}\overset{}{_{fd}^{a}} & - & \underset{\left(
11\right) }{\overset{V}{C}}\overset{}{_{bd}^{f}}\underset{\left( 11\right) }{%
\overset{V}{C}}\overset{}{_{fc}^{a}}. & 
\end{array}
\label{curb_var_v}
\end{equation}%

\section{The relative covariant derivatives}

Let $D\Gamma \left( N\right) $, the Cartan metrical N-linear connection of the manifold $OscM$. A
classical method to determine the laws of derivation on a Finsler
submanifold is the type of the coupling$.$ \newline
\textbf{Theorem 2.1 }\textit{The coupling of the N-linear connection D\ to
the induced nonli\-near connection \v{N} along }$\widetilde{Osc\check{M}}$\textit{\ is
locally given by\ the set of coefficients} $\check{D}\Gamma \left( \check{N}%
\right) =\left( \underset{\left( 00\right) }{\overset{H}{\check{L}}}\overset{%
}{_{b\delta }^{a}},\underset{\left( 10\right) }{\overset{V}{\check{L}}}%
\overset{}{_{b\delta }^{a}},\underset{\left( 01\right) }{\overset{H}{\check{C%
}}}\overset{}{_{b\delta }^{a}},\underset{\left( 11\right) }{\overset{V}{%
\check{C}}}\overset{}{_{b\delta }^{a}}\right) ,$ \textit{where}%
\begin{equation}
\left\{ 
\begin{array}{l}
\underset{\left( 00\right) }{\overset{H}{\check{L}}}\overset{}{_{b\delta
}^{a}}=\underset{\left( 00\right) }{\overset{H}{L}}\overset{}{_{bd}^{a}}%
B_{\delta }^{d}+\underset{\left( 01\right) }{\overset{H}{C}}\overset{}{%
_{bd}^{a}}B_{\bar{\delta}}^{d}K\overset{}{_{\delta }^{\bar{\delta}}} \\ 
\\ 
\underset{\left( 10\right) }{\overset{V}{\check{L}}}\overset{}{_{b\delta
}^{a}}=\underset{\left( 10\right) }{\overset{V}{L}}\overset{}{_{bd}^{a}}%
B_{\delta }^{d}+\underset{\left( 11\right) }{\overset{V}{C}}\overset{}{%
_{bd}^{a}}B_{\bar{\delta}}^{d}K\overset{}{_{\delta }^{\bar{\delta}}} \\ 
\\ 
\underset{\left( 01\right) }{\overset{H}{\check{C}}}\overset{}{_{b\delta
}^{a}}=\underset{\left( 01\right) }{\overset{H}{C}}\overset{}{_{bd}^{a}}%
B_{\delta }^{d} \\ 
\\ 
\underset{\left( 11\right) }{\overset{V}{\check{C}}}\overset{}{_{b\delta
}^{a}}=\underset{\left( 11\right) }{\overset{V}{C}}\overset{}{_{bd}^{a}}%
B_{\delta }^{d}.%
\end{array}%
\right.  \label{coef_rac}
\end{equation}%
\textbf{Definition 2.2 }\ \textit{We call the \textbf{induced tangent
connection} on }$\widetilde{Osc\check{M}}$\textit{\ by the metrical N-linear connection
D, the couple of the operators }$\overset{H}{D^{\top }}$, $\overset{V}{%
D^{\top }}$ \textit{which are defined by}%
\begin{equation*}
\begin{array}{c}
\begin{array}{l}
\overset{H}{D^{\top }}X^{\alpha }=B_{b}^{\alpha }\overset{H}{\check{D}}X^{b},
\\ 
\\ 
\overset{V}{D^{\top }}X^{\alpha }=B_{b}^{\alpha }\overset{V}{\check{D}}X^{b},%
\end{array}%
\quad \textit{for }X\overset{}{^{a}}=B_{\gamma }^{a}X^{\gamma }%
\end{array}%
\end{equation*}
\textit{where} 
\begin{equation*}
\begin{array}{l}
\overset{H}{D^{\top }}X^{\alpha }=dX^{\alpha }+X^{\beta }\overset{H}{\omega }%
\overset{}{_{\beta }^{\alpha }} \\ 
\\ 
\overset{V}{D^{\top }}X^{\alpha }=dX^{\alpha }+X^{\beta }\overset{V}{\omega }%
\overset{}{_{\beta }^{\alpha }}%
\end{array}%
\end{equation*}
\textit{and }$\overset{H}{\omega }\overset{}{_{\beta }^{\alpha }}$\textit{, }%
$\overset{V}{\omega }\overset{}{_{\beta }^{\alpha }}$\textit{\ are called
the }\textbf{\textit{tangent connection 1-forms}}$.$

We have \newline
\textbf{Theorem 2.3} \ \textit{The tangent connections 1-forms are as
follows:}%
\begin{equation}
\begin{array}{l}
\overset{H}{\omega }\overset{}{_{\beta }^{\alpha }}=\underset{\left(
00\right) }{\overset{H}{L}}\overset{}{_{\beta \delta }^{\alpha }}du^{\delta
}+\underset{\left( 01\right) }{\overset{H}{C}}\overset{}{_{\beta \delta
}^{\alpha }}\delta v^{\delta } \\ 
\\ 
\overset{V}{\omega }\overset{}{_{\beta }^{\alpha }}=\underset{\left(
10\right) }{\overset{V}{L}}\overset{}{_{\beta \delta }^{\alpha }}du^{\delta
}+\underset{\left( 11\right) }{\overset{V}{C}}\overset{}{_{\beta \delta
}^{\alpha }}\delta v^{\delta },%
\end{array}
\label{1-form_con_tg}
\end{equation}%
\textit{where}%
\begin{equation}
\begin{array}{l}
\underset{\left( 00\right) }{\overset{H}{L}}\overset{}{_{\beta \delta
}^{\alpha }}=B_{d}^{\alpha }\left( B_{\beta \delta }^{d}+B_{\beta }^{f}%
\underset{\left( 00\right) }{\overset{H}{\check{L}}}\overset{}{_{f\delta
}^{d}}\right) , \\ 
\\ 
\underset{\left( 10\right) }{\overset{V}{L}}\overset{}{_{\beta \delta
}^{\alpha }}=B_{d}^{\alpha }\left( B_{\beta \delta }^{d}+B_{\beta }^{f}%
\underset{\left( 10\right) }{\overset{V}{\check{L}}}\overset{}{_{f\delta
}^{d}}\right) , \\ 
\\ 
\underset{\left( 01\right) }{\overset{H}{C}}\overset{}{_{\beta \delta
}^{\alpha }}=B_{d}^{\alpha }B_{\beta }^{f}\underset{\left( 01\right) }{%
\overset{H}{\check{C}}}\overset{}{_{f\delta }^{d}}, \\ 
\\ 
\underset{\left( 11\right) }{\overset{V}{C}}\overset{}{_{\beta \delta
}^{\alpha }}=B_{d}^{\alpha }B_{\beta }^{f}\underset{\left( 11\right) }{%
\overset{V}{\check{C}}}\overset{}{_{f\delta }^{d}}.%
\end{array}
\label{coef_con_tg}
\end{equation}%
\textbf{Definition 2.4} \textit{We call the \textbf{induced normal connection%
} on }$\widetilde{Osc\check{M}}$\textit{\ by the metrical N-linear connection D, the
couple of the operators }$\overset{H}{D^{\bot }}$, $\overset{V}{D^{\bot }}$ 
\textit{which are defined by}%
\begin{equation*}
\begin{array}{c}
\begin{array}{l}
\overset{H}{D^{\bot }}X^{\overline{\alpha }}=B_{b}^{\alpha }\overset{H}{%
\check{D}}X^{b} \\ 
\\ 
\overset{V}{D^{\bot }}X^{\overline{\alpha }}=B_{b}^{\alpha }\overset{V}{%
\check{D}}X^{b},%
\end{array}%
\quad \text{\textit{for} }X^{a}=B_{\bar{\gamma}}^{a}X^{\bar{\gamma}}%
\end{array}%
\end{equation*}
\textit{where} 
\begin{equation*}
\begin{array}{l}
\overset{H}{D^{\bot }}X^{\overline{\alpha }}=dX^{\overline{\alpha }}+X^{%
\overline{\beta }}\overset{H}{\omega }\overset{}{_{\overline{\beta }}^{%
\overline{\alpha }}} \\ 
\\ 
\overset{V}{D^{\bot }}X^{\overline{\alpha }}=dX^{\overline{\alpha }}+X^{%
\overline{\beta }}\overset{V}{\omega }\overset{}{_{\overline{\beta }}^{%
\overline{\alpha }}}%
\end{array}%
\end{equation*}
\textit{and }$\overset{H}{\omega }\overset{}{_{\overline{\beta }}^{\overline{%
\alpha }},\text{ }\overset{V}{\omega }\overset{}{_{\overline{\beta }}^{%
\overline{\alpha }}}}$\textit{\ are called the }\textbf{\textit{normal
connection 1-forms}}$.$

We have \newline
\textbf{Theorem 2.5} \ \textit{The normal connections 1-forms are as follows:%
}

\begin{equation}
\begin{array}{l}
\overset{H}{\omega }\overset{}{_{\overline{\beta }}^{\overline{\alpha }}}=%
\underset{\left( 00\right) }{\overset{H}{L}}\overset{}{_{\overline{\beta }%
\delta }^{\overline{\alpha }}}du^{\delta }+\underset{\left( 01\right) }{%
\overset{H}{C}}\overset{}{_{\overline{\beta }\delta }^{\overline{\alpha }}}%
\delta v^{\delta } \\ 
\\ 
\overset{V}{\omega }\overset{}{_{\overline{\beta }}^{\overline{\alpha }}}=%
\underset{\left( 10\right) }{\overset{V}{L}}\overset{}{_{\overline{\beta }%
\delta }^{\overline{\alpha }}}du^{\delta }+\underset{\left( 11\right) }{%
\overset{V}{C}}\overset{}{_{\overline{\beta }\delta }^{\overline{\alpha }}}%
\delta v^{\delta },%
\end{array}
\label{1-form_con_norm}
\end{equation}%
\newline
\textit{where}

\begin{equation}
\begin{array}{l}
\underset{\left( 00\right) }{\overset{H}{L}}\overset{}{_{\overline{\beta }%
\delta }^{\overline{\alpha }}}=B_{d}^{\overline{\alpha }}\left( \dfrac{%
\delta B_{\overline{\beta }}^{d}}{\delta u^{\delta }}+B_{\overline{\beta }%
}^{f}\underset{\left( 00\right) }{\overset{H}{\check{L}}}\overset{}{%
_{f\delta }^{d}}\right) , \\ 
\\ 
\underset{\left( 10\right) }{\overset{V}{L}}\overset{}{_{\overline{\beta }%
\delta }^{\overline{\alpha }}}=B_{d}^{\overline{\alpha }}\left( \dfrac{%
\delta B_{\overline{\beta }}^{d}}{\delta u^{\delta }}+B_{\overline{\beta }%
}^{f}\underset{\left( 10\right) }{\overset{V}{\check{L}}}\overset{}{%
_{f\delta }^{d}}\right) , \\ 
\\ 
\underset{\left( 01\right) }{\overset{H}{C}}\overset{}{_{\overline{\beta }%
\delta }^{\overline{\alpha }}}=B_{d}^{\overline{\alpha }}\left( \dfrac{%
\partial B_{\overline{\beta }}^{d}}{\partial v^{\delta }}+B_{\overline{\beta 
}}^{f}\underset{\left( 01\right) }{\overset{H}{\check{C}}}\overset{}{%
_{f\delta }^{d}}\right) , \\ 
\\ 
\underset{\left( 11\right) }{\overset{V}{C}}\overset{}{_{\overline{\beta }%
\delta }^{\overline{\alpha }}}=B_{d}^{\overline{\alpha }}\left( \dfrac{%
\partial B_{\overline{\beta }}^{d}}{\partial v^{\delta }}+B_{\overline{\beta 
}}^{f}\underset{\left( 11\right) }{\overset{V}{\check{C}}}\overset{}{%
_{f\delta }^{d}}\right) .%
\end{array}
\label{coef_con_norm}
\end{equation}

Now, we can define the relative (or mixed) covariant derivatives $\overset{H}{\nabla }
$ and $\overset{V}{\nabla }$ . 
\newline
\textbf{Theorem 2.6}\ \textit{The relative covariant (mixed) derivatives in
the algebra of mixed d-tensor fields are the operators }$\overset{H}{\nabla }
$, $\overset{V}{\nabla }$\textit{\ for which the following properties hold}: 
\begin{equation*}
\begin{array}{c}
\begin{array}{l}
\overset{H}{\nabla }f=df, \\ 
\\ 
\overset{V}{\nabla }f=df,%
\end{array}%
\quad \forall f\in \mathcal{F}\left( \widetilde{Osc\check{M}}\right) \\ 
\\ 
\begin{array}{l}
\overset{H}{\nabla }X^{a}=\overset{H}{\check{D}}X^{a}, \\ 
\\ 
\overset{V}{\nabla }X^{a}=\overset{V}{\check{D}}X^{a},%
\end{array}%
\begin{array}{l}
\overset{H}{\nabla }X^{\alpha }=\overset{H}{D^{\intercal }}X^{\alpha }, \\ 
\\ 
\overset{V}{\nabla }X^{\alpha }=\overset{V}{D^{\intercal }}X^{\alpha },%
\end{array}%
\begin{array}{l}
\overset{H}{\nabla }X^{\overline{\alpha }}=\overset{H}{D^{\bot }}X^{%
\overline{\alpha }}, \\ 
\\ 
\overset{V}{\nabla }X^{\overline{\alpha }}=\overset{H}{D^{\bot }}X^{%
\overline{\alpha }}.%
\end{array}%
\end{array}%
\end{equation*}

$\overset{H}{\check{\omega}}\overset{}{_{b}^{a}},\overset{V}{\check{\omega}}%
\overset{}{_{b}^{a}},\overset{H}{\omega }\overset{}{_{\beta }^{\alpha }},%
\overset{V}{\omega }\overset{}{_{\beta }^{\alpha }},\overset{H}{\omega }%
\overset{}{_{\overline{\beta }}^{\overline{\alpha }}}$, $\overset{V}{\omega }%
\overset{}{_{\overline{\beta }}^{\overline{\alpha }}}$ \textit{are called the%
} \textit{\textbf{connection 1-forms} of }$\overset{H}{\nabla }$, $\overset{V%
}{\nabla }$\textit{$.$}

\section{A comparison between the induced and intrinsic geometrical objects}

It is known that, in the case of Finsler or pseudo-Finsler spaces (\cite{Ma}%
, \cite{Bej+Far}, \cite{Munt1}, \cite{Munt2}, \cite{Munt3}, \cite{Munt4}),\
the intrinsic nonlinear connection of the submanifold is \linebreak
different from the induced nonlinear connection by the nonlinear connection
on the manifold. Moreover, the induced Finsler connection is different from
the the induced Finsler connection.

In this section, we present a comparison between the induced and intrinsic
geometric objects on the submanifold $\widetilde{Osc}\check{M}$ with respect
to the Cartan metrical $N$-linear connection determinated by the homogeneous
lift $\overset{0}{\mathbb{G}}$ (\ref{lift_om}).

Let $\widetilde{\overset{0}{\mathbb{G}}},$ the homogeneous lift to the
submanifold $\widetilde{Osc\check{M}}$ of the induced fundamental tensor (%
\ref{tens_fund_subvar}), 
\begin{equation}
\begin{array}{r}
\widetilde{\overset{0}{\mathbb{G}}}=g_{\alpha \beta }\left( u,v\right)
du^{\alpha }\otimes du^{\beta }+h_{\alpha \beta }\left( u,v\right) \delta
v^{\alpha }\otimes \delta v^{\beta },%
\end{array}
\label{lift_indus}
\end{equation}%
where 
\begin{equation*}
\begin{array}{l}
h_{\alpha \beta }\left( u,v\right) =\dfrac{p^{2}}{\left\Vert v\right\Vert
^{2}}g_{\alpha \beta }\left( u,v\right) \\ 
\\ 
\left\Vert v\right\Vert ^{2}=g_{\alpha \beta }v^{\alpha }v^{\beta }%
\end{array}%
\end{equation*}%
and $\mathring{N}$, the intrinsic Cartan nonlinear connection 
\begin{equation*}
\mathring{N}\overset{}{^{\alpha }}_{\beta }=\dfrac{\partial G^{\alpha }}{%
\partial v^{\beta }},
\end{equation*}%
where $G^{\alpha }=\dfrac{1}{2}\gamma \overset{}{_{\beta \gamma }^{\alpha }}%
\left( u,v\right) v^{\beta }v^{\gamma }$ and $\gamma \overset{}{_{\beta
\gamma }^{\alpha }}\left( u,v\right) $ are the Christoffel symbols of $%
g_{\alpha \beta }.$

Let $\mathring{D}\Gamma \left( \mathring{N}\right) =\left( \underset{\left(
00\right) }{\overset{H}{\mathring{L}}}\overset{}{_{\beta \gamma }^{\alpha }},%
\underset{\left( 10\right) }{\overset{V_{1}}{\mathring{L}}}\overset{}{%
_{\beta \gamma }^{\alpha }},\underset{\left( 01\right) }{\overset{H}{%
\mathring{C}}}\overset{}{_{\beta \gamma }^{\alpha }},\underset{\left(
11\right) }{\overset{V_{1}}{\mathring{C}}}\overset{}{_{\beta \gamma
}^{\alpha }}\right) $ the intrinsic Cartan metrical \linebreak $N$-linear
connection of the submanifold $\widetilde{Osc\check{M}}$ with \newline
the "\textit{horizontal}" coeficients

\begin{equation}
\begin{array}{l}
\underset{\left( 00\right) }{\overset{H}{\mathring{L}}}\overset{}{_{\beta
\gamma }^{\alpha }}=\dfrac{1}{2}g^{\alpha \delta }\left( \delta _{\beta
}g_{\delta \gamma }+\delta _{\gamma }g_{\beta \delta }-\delta _{\delta
}g_{\beta \gamma }\right) \\ 
\\ 
\underset{\left( 10\right) }{\overset{V}{\mathring{L}}}\overset{}{_{\beta
\gamma }^{\alpha }}=\dfrac{1}{2}h^{\alpha \delta }\left( \delta _{\beta
}h_{\delta \gamma }+\delta _{\gamma }h_{\beta \delta }-\delta _{\delta
}h_{\beta \gamma }\right)%
\end{array}
\label{coef_N_con_subvar_h}
\end{equation}%
and the "\textit{vertical}" coeficients

\begin{equation}
\begin{array}{l}
\underset{\left( 01\right) }{\overset{H}{\mathring{C}}}\overset{}{_{\beta
\gamma }^{\alpha }}=\dfrac{1}{2}g^{\alpha \delta }\left( \dot{\partial}%
_{1\beta }g_{\delta \gamma }+\dot{\partial}_{1\gamma }g_{\beta \delta }-\dot{%
\partial}_{1\delta }g_{\beta \gamma }\right) \\ 
\\ 
\underset{\left( 11\right) }{\overset{V}{\mathring{C}}}\overset{}{_{\beta
\gamma }^{\alpha }}=\dfrac{1}{2}h^{\alpha \delta }\left( \dot{\partial}%
_{1\beta }h_{\delta \gamma }+\dot{\partial}_{1\gamma }h_{\beta \delta }-\dot{%
\partial}_{1\delta }h_{\beta \gamma }\right) .%
\end{array}
\label{coef_N_con_subvar_v}
\end{equation}%
\newline
\textbf{Proposition 3.1} \textit{The Lie brakets of the vector fields of the
adapted basis }$\left\{ \delta _{\alpha },\dot{\partial}_{1\alpha }\right\} $%
\textit{\ are given by}:%
\begin{equation}
\begin{array}{l}
\left[ \delta _{\alpha },\delta _{\beta }\right] =-R\underset{}{_{\alpha
\beta }^{\sigma }}\dot{\partial}_{1\sigma }, \\ 
\\ 
\left[ \delta _{\alpha },\dot{\partial}_{1\beta }\right] =\underset{\left(
11\right) }{B}\underset{}{_{\alpha \beta }^{\sigma }}\dot{\partial}_{1\sigma
}, \\ 
\\ 
\left[ \dot{\partial}_{1\alpha },\dot{\partial}_{1\beta }\right] =0,%
\end{array}%
\end{equation}%
\label{Lie_subvar} \textit{where}%
\begin{equation*}
\begin{array}{lll}
R\underset{}{_{\alpha \beta }^{\sigma }} & = & \delta _{\alpha }N\overset{}{%
^{\sigma }}_{\beta }-\delta _{\beta }N\overset{}{^{\sigma }}_{\alpha } \\ 
&  &  \\ 
\underset{\left( 11\right) }{B}\underset{}{_{\alpha \beta }^{\sigma }} & = & 
\dot{\partial}_{1\beta }N\overset{}{^{\sigma }}_{\alpha }.%
\end{array}%
\end{equation*}

For any d-vector field $X\in X(\widetilde{OscM})$ expressed in the adapted basis $%
\left\{ \delta _{a},\dot{\partial}_{1a}\right\} $ we have%
\begin{equation*}
X=\overset{0}{X}\underset{}{^{a}}\dfrac{\delta }{\delta x^{a}}+\overset{1}{X}%
\underset{}{^{a}}\dfrac{\partial }{\partial y^{a}},X\in \mathcal{X}\left( 
\widetilde{OscM}\right) .
\end{equation*}%
We consider $h$ and $v$, the horizontal and the vertical projectors
associated to the nonlinear connection $N.$ Denote by

\begin{eqnarray*}
X^{H} &=&hX=\overset{0}{X}\underset{}{^{a}}\dfrac{\delta }{\delta x^{a}} \\
&& \\
X^{V} &=&vX=\overset{1}{X}\underset{}{^{a}}\dfrac{\partial }{\partial y^{a}}.
\end{eqnarray*}%
For any d-vector field $\check{X}\in X(\widetilde{Osc\check{M}})$ expressed in the
adapted basis $\left\{ \delta _{\alpha },\dot{\partial}_{1\alpha }\right\} $
we have

\begin{equation*}
\check{X}=\overset{0}{\check{X}}\underset{}{^{\alpha }}\dfrac{\mathring{%
\delta}}{\mathring{\delta}u^{\alpha }}+\overset{1}{\check{X}}\underset{}{%
^{\alpha }}\dfrac{\partial }{\partial v^{\alpha }},\forall \check{X}\in 
\mathcal{X}\left( \widetilde{Osc\check{M}}\right)
\end{equation*}%
and consider $\mathring{h}$ and $\mathring{v},$ the horizontal and the
vertical projectors associated to the intrinsic nonlinear connection $\mathring{N}$.
Denote by%
\begin{equation*}
\begin{array}{c}
\check{X}^{\mathring{H}}=\mathring{h}\check{X}=\overset{0}{\check{X}}%
\underset{}{^{\alpha }}\dfrac{\mathring{\delta}}{\mathring{\delta}u^{\alpha }%
} \\ 
\\ 
\check{X}^{\mathring{V}}=v\check{X}=\overset{1}{\check{X}}\underset{}{^{a}}%
\dfrac{\partial }{\partial v^{a}}%
\end{array}%
\forall \check{X}\in \mathcal{X}\left( \widetilde{Osc\check{M}}\right) .
\end{equation*}
\newline
\textbf{Proposition 3.2} \textit{Let }$\mathring{N}$\textit{, the intrinsic
Cartan nonlinear connection and }$\check{N}$\textit{, the induced nonlinear
connection on the submanifold }$\widetilde{Osc\check{M}}$\textit{\ by the
Cartan nonlinear connection N. Then the following relations hold}:

\bigskip \textit{1}$^{\circ }$\textit{The coefficients of the nonlinear
connections }$\mathring{N}$\textit{\ and }$\check{N}$\textit{\ are related
by (\cite{Bej+Far})}%
\begin{equation*}
\mathring{N}^{\alpha }\overset{}{_{\beta }}=\check{N}^{\alpha }\overset{}{%
_{\beta }}+D^{\alpha }\overset{}{_{\beta }}.
\end{equation*}

\textit{2}$^{\circ }$\textit{There exist the following relations between the
components of the adapted bases of }$\mathring{N}$\textit{\ and }$\check{N}$%
\begin{equation*}
\begin{array}{l}
\mathring{\delta}_{\alpha }=\delta _{\alpha }-\overset{}{D^{\beta }}\overset{%
}{_{\alpha }}\dot{\partial}_{1\beta }, \\ 
\\ 
\mathring{\partial}_{1\alpha }=\dot{\partial}_{1\alpha }.%
\end{array}%
\end{equation*}

\textit{3}$^{\circ }$\textit{There exist the following relations between the
coefficients of the Lie brakets of the adapted bases of \r{N} and \v{N}}:%
\begin{equation*}
\begin{array}{l}
\mathring{R}\overset{}{_{\beta \gamma }^{\alpha }}=R\overset{}{_{\beta
\gamma }^{\alpha }}+\overset{}{\underset{00}{\overset{1}{D}}}\underset{}{%
_{\beta \gamma }^{\alpha }} \\ 
\\ 
\underset{\left( 11\right) }{\mathring{B}}\overset{}{_{\beta \gamma
}^{\alpha }}=\underset{\left( 11\right) }{B}\overset{}{_{\beta \gamma
}^{\alpha }}+\overset{}{\underset{01}{\overset{1}{D}}}\underset{}{_{\beta
\gamma }^{\alpha }},%
\end{array}%
\end{equation*}%
\textit{where}%
\begin{equation}
\begin{array}{lll}
D\overset{}{^{\alpha }}\overset{}{_{\beta }} & = & g\underset{}{_{\bar{\alpha%
}}}\underset{}{^{\alpha }}\underset{}{_{\beta }}K\underset{}{_{\beta }^{\bar{%
\alpha}}}v^{\beta }%
\end{array}
\label{rest_con_nel}
\end{equation}

\begin{equation}
\overset{}{\underset{00}{\overset{1}{D}}}\underset{}{_{\beta \gamma
}^{\alpha }}=\delta _{\gamma }D^{\alpha }\overset{}{_{\beta }}-\delta
_{\beta }D^{\alpha }\overset{}{_{\gamma }}+D^{\delta }\overset{}{_{\beta }}%
\dot{\partial}_{1\delta }\left( N^{\alpha }\overset{}{_{\gamma }}+D^{\alpha }%
\overset{}{_{\gamma }}\right) -D^{\delta }\overset{}{_{\gamma }}\dot{\partial%
}_{1\delta }\left( N^{\alpha }\overset{}{_{\beta }}+D^{\alpha }\overset{}{%
_{\beta }}\right)  \label{D100}
\end{equation}

\begin{equation}
\overset{}{\overset{}{\underset{01}{\overset{1}{D}}}\underset{}{_{\beta
\gamma }^{\alpha }}=}\dot{\partial}_{1\gamma }D^{\alpha }\overset{}{_{\beta }%
}.  \label{D101}
\end{equation}%
\newline
\textbf{Proposition 3.3} \textit{The local coefficients of the intrinsic
Cartan metrical }$N$\textit{-linear connection }$\mathring{D}\Gamma \left( 
\mathring{N}\right) $\textit{\ and of the induced tangent connection of the
Cartan \linebreak metrical }$N$\textit{-linear connection }$D\Gamma \left(
N\right) $ \textit{are related by}:%
\begin{equation*}
\begin{array}{l}
\underset{\left( 00\right) }{\overset{H}{\mathring{L}}}\overset{}{_{\beta
\delta }^{\alpha }}=\underset{\left( 00\right) }{\overset{H}{L}}\underset{}{%
_{\beta \delta }^{\alpha }}+\underset{\left( 00\right) }{\overset{H}{\Delta }%
}\underset{}{_{\beta \delta }^{\alpha }} \\ 
\\ 
\underset{\left( 10\right) }{\overset{V}{\mathring{L}}}\overset{}{_{\beta
\delta }^{\alpha }}=\underset{\left( 10\right) }{\overset{V}{L}}\overset{}{%
_{\beta \delta }^{\alpha }}+\underset{\left( 10\right) }{\overset{V}{\Delta }%
}\underset{}{_{\beta \delta }^{\alpha }}%
\end{array}%
\end{equation*}

\begin{equation*}
\begin{array}{l}
\underset{\left( 01\right) }{\overset{H}{\mathring{C}}}\overset{}{_{\beta
\delta }^{\alpha }}=\underset{\left( 01\right) }{\overset{H}{C}}\overset{}{%
_{\beta \delta }^{\alpha }}, \\ 
\\ 
\underset{\left( 11\right) }{\overset{V}{\mathring{C}}}\overset{}{_{\beta
\delta }^{\alpha }}=\underset{\left( 11\right) }{\overset{V}{C}}\overset{}{%
_{\beta \delta }^{\alpha }},%
\end{array}%
\end{equation*}
\textit{where}%
\begin{equation}
\begin{array}{ccl}
\underset{\left( 00\right) }{\overset{H}{\Delta }}\underset{}{_{\beta \delta
}^{\alpha }} & = & \dfrac{1}{2}\left[ g\underset{}{_{\bar{\delta}}}\underset{%
}{^{\alpha }}\underset{}{_{\delta }}K\underset{}{_{\beta }^{\bar{\delta}}}-g%
\underset{}{_{\bar{\delta}\beta \delta }}K\underset{}{_{{}}^{\bar{\delta}%
\alpha }}\right] \\ 
&  &  \\ 
&  & -\dfrac{1}{2}g^{\alpha \delta }\left[ D\overset{}{^{\varepsilon }}%
\overset{}{_{\beta }}\delta _{1\varepsilon }g_{\delta \gamma }+D\overset{}{%
^{\varepsilon }}\overset{}{_{\gamma }}\delta _{1\varepsilon }g_{\delta \beta
}-D\overset{}{^{\varepsilon }}\overset{}{_{\delta }}\delta _{1\varepsilon
}g_{\beta \gamma }\right] , \\ 
&  &  \\ 
\underset{\left( 10\right) }{\overset{V}{\Delta }}\underset{}{_{\beta \delta
}^{\alpha }} & = & -\dfrac{1}{2h}g^{\alpha \sigma }\left[ B_{\bar{\beta}%
}^{b}K_{\sigma }^{\bar{\beta}}\left( \dot{\partial}_{1b}h\right) g_{\beta
\delta }-B_{\bar{\beta}}^{b}K_{\beta }^{\bar{\beta}}\left( \dot{\partial}%
_{1b}h\right) g_{\sigma \delta }\right] \\ 
&  &  \\ 
&  & +\dfrac{1}{2h}g^{\alpha \sigma }\left[ D^{\varepsilon }\overset{}{%
_{\sigma }}\dot{\partial}_{1\varepsilon }g_{\beta \delta }-D^{\varepsilon }%
\overset{}{_{\beta }}\dot{\partial}_{1\varepsilon }g_{\sigma \delta
}-D^{\varepsilon }\overset{}{_{\delta }}\dot{\partial}_{1\varepsilon
}g_{\beta \sigma }\right] \\ 
&  &  \\ 
&  & +\underset{\left( 10\right) }{\Delta }\underset{}{_{\beta \delta
}^{\alpha }},%
\end{array}
\label{deltaHV}
\end{equation}
\textit{and}%
\begin{equation}
\begin{array}{ccl}
\underset{\left( 10\right) }{\Delta }\underset{}{_{\beta \gamma }^{\alpha }}
& = & \dfrac{1}{2}\left( \dot{\partial}_{1\beta }B\underset{}{_{a}^{\alpha }}%
\right) \left( B\underset{}{_{0\gamma }^{a}}+B\underset{}{_{\gamma }^{c}}%
\underset{1}{N}\underset{}{^{a}}\underset{}{_{c}}\right) -\dfrac{1}{2}B%
\underset{}{_{a}^{\alpha }}B\underset{}{_{\beta \gamma }^{a}} \\ 
&  &  \\ 
&  & -\dfrac{1}{2}g^{af}B\underset{}{_{a}^{\alpha }}B\underset{}{_{\beta
}^{b}}B\underset{}{_{\overline{\delta }}^{d}}\left( \dot{\partial}%
_{1d}g_{bf}\right) K\underset{}{^{\overline{\delta }}}\underset{}{_{\gamma }}%
-\dfrac{1}{2}g^{\alpha \delta }B\underset{}{_{a}^{\sigma }}B\underset{}{%
_{\gamma \delta }^{a}}g_{\sigma \beta } \\ 
&  &  \\ 
&  & -\dfrac{1}{2}g^{\alpha \delta }g_{\sigma \beta }\left( \dot{\partial}%
_{1\gamma }B\underset{}{_{a}^{\sigma }}\right) \left( B\underset{}{_{0\delta
}^{a}}+B\underset{}{_{\delta }^{d}}N\underset{}{^{a}}\underset{}{_{d}}\right)
\\ 
&  &  \\ 
&  & +g^{\alpha \delta }g_{bd}B\underset{}{_{\beta }^{b}}B\underset{}{%
_{\delta \gamma }^{a}}+\dfrac{1}{2}\underset{01}{\overset{1}{D}}\underset{}{%
_{\gamma \beta }^{\alpha }}-\dfrac{1}{2}g^{\alpha \delta }\left[ D\underset{}%
{^{\varepsilon }}\underset{}{_{\gamma }}\dot{\partial}_{1\varepsilon
}g_{\beta \delta }+\underset{01}{\overset{1}{D}}\underset{}{_{\gamma \delta
}^{\sigma }}g_{\sigma \beta }\right] .%
\end{array}
\label{delta10}
\end{equation}

From the proposition 3.3 we get that $\mathring{D}$, the intrinsic Cartan
metrical \linebreak N-linear connection is not identical with $D^{\top },$
the induced tangent connection of the Cartan metrical N-linear connection D$%
\Gamma \left( N\right) .$ From this fact, \cite{Bej+Far}, there exist $\overset{\top }{D}$,
the \textbf{deformation tensor} of the pair $(\mathring{D},D^{\top })$. 
For $\check{X},\check{Y}\in \mathcal{\check{X}}\left( \widetilde{Osc\check{M}}%
\right) $ we get

\begin{equation*}
\begin{array}{l}
\begin{array}{c}
\mathring{D}_{\check{X}^{\mathring{V}}}\check{Y}^{\mathring{H}}=D_{\check{X}%
^{\mathring{V}}}^{\top }\check{Y}^{\mathring{H}} \\ 
\\ 
\mathring{D}_{\check{X}^{\mathring{V}}}\check{Y}^{\mathring{V}}=D_{\check{X}%
^{\mathring{V}}}^{\top }\check{Y}^{\mathring{V}}%
\end{array}
\\ 
\\ 
\begin{array}{c}
\mathring{D}_{\check{X}^{\mathring{H}}}\check{Y}^{\mathring{H}}=D_{\check{X}%
^{\mathring{H}}}^{\top }\check{Y}^{\mathring{H}}+\overset{\top }{D}\left( 
\check{X}^{\mathring{H}},\check{Y}^{\mathring{H}}\right) \\ 
\\ 
\mathring{D}_{\check{X}^{\mathring{H}}}\check{Y}^{\mathring{V}}=D_{\check{X}%
^{\mathring{H}}}^{\top }\check{Y}^{\mathring{V}}+\overset{\top }{D}\left( 
\check{X}^{\mathring{H}},\check{Y}^{\mathring{V}}\right) ,%
\end{array}%
\end{array}%
\begin{array}{c}
\\ 
\forall \check{X},\check{Y}\in \mathcal{\check{X}}\left( \widetilde{Osc%
\check{M}}\right) . \\ 
\end{array}%
\end{equation*}%

If we express $\overset{\top }{D}$
in the adapted bases of $\mathring{N},$ we get:%
\begin{equation*}
\begin{array}{c}
\overset{\top }{D}\left( \dfrac{\mathring{\delta}}{\mathring{\delta}%
u^{\gamma }},\dfrac{\mathring{\delta}}{\mathring{\delta}u^{\beta }}\right) =%
\underset{00}{\overset{\top }{D}}\underset{}{^{\mathring{H}}}\overset{}{%
_{\beta \gamma }^{\alpha }}\dfrac{\mathring{\delta}}{\mathring{\delta}%
u^{\alpha }}+\underset{00}{\overset{\top }{D}}\underset{}{^{\mathring{V}}}%
\overset{}{_{\beta \gamma }^{\alpha }}\dfrac{\partial }{\partial v^{a}} \\ 
\\ 
\overset{\top }{D}\left( \dfrac{\mathring{\delta}}{\mathring{\delta}%
u^{\gamma }},\dfrac{\partial }{\partial v^{\beta }}\right) =\underset{10}{%
\overset{\top }{D}}\underset{}{^{\mathring{H}}}\overset{}{_{\beta \gamma
}^{\alpha }}\dfrac{\mathring{\delta}}{\mathring{\delta}u^{\alpha }}+\underset%
{10}{\overset{\top }{D}}\underset{}{^{\mathring{V}}}\overset{}{_{\beta
\gamma }^{\alpha }}\dfrac{\partial }{\partial v^{a}}.%
\end{array}%
\end{equation*}

We have the next \newline
\textbf{Proposition 3.4} \textit{The components of the deformation tensor }$%
\overset{\top }{D}$\textit{\ are given by the formula}:%
\begin{equation*}
\begin{array}{ccl}
\underset{00}{\overset{\top }{D}}\underset{}{^{\mathring{H}}}\overset{}{%
_{\beta \gamma }^{\alpha }} & = & \underset{\left( 00\right) }{\overset{H}{%
\Delta }}\overset{}{_{\beta \gamma }^{\alpha }}+D_{\gamma }^{\varphi }%
\underset{\left( 01\right) }{\overset{H}{C}}\overset{}{_{\beta \varphi
}^{\alpha }} \\ 
&  &  \\ 
\underset{00}{\overset{\top }{D}}\underset{}{^{\mathring{V}}}\overset{}{%
_{\beta \gamma }^{\varepsilon }} & = & \underset{00}{\overset{\top }{D}}%
\underset{}{^{\mathring{H}}}\overset{}{_{\beta \gamma }^{\varepsilon }}%
D_{\varepsilon }^{\alpha }+\delta _{\gamma }D_{\beta }^{\alpha }+D_{\beta
}^{\varphi }\underset{\left( 10\right) }{\overset{V}{L}}\overset{}{_{\varphi
\gamma }^{\alpha }}-D_{\gamma }^{\varepsilon }\left[ \dot{\partial}%
_{1\varepsilon }D_{\beta }^{\alpha }+D_{\beta }^{\varphi }\underset{\left(
11\right) }{\overset{V}{C}}\overset{}{_{\varphi \varepsilon }^{\alpha }}%
\right] \\ 
&  &  \\ 
&  & -D_{\alpha }^{\varepsilon }\left( \underset{\left( 00\right) }{\overset{%
H}{L}}\overset{}{_{\beta \gamma }^{\alpha }}+\underset{\left( 00\right) }{%
\Delta }\overset{}{_{\beta \gamma }^{\alpha }}\right) , \\ 
&  &  \\ 
\underset{10}{\overset{\top }{D}}\underset{}{^{\mathring{H}}}\overset{}{%
_{\beta \gamma }^{\alpha }} & = & 0, \\ 
&  &  \\ 
\underset{10}{\overset{\top }{D}}\underset{}{^{\mathring{V}}}\overset{}{%
_{\beta \gamma }^{\alpha }} & = & \underset{\left( 10\right) }{\overset{V}{%
\Delta }}\overset{}{_{\beta \gamma }^{\alpha }}+D_{\gamma }^{\varphi }%
\underset{\left( 01\right) }{\overset{H}{C}}\overset{}{_{\beta \varphi
}^{\alpha }},%
\end{array}%
\end{equation*}%
\textit{where} $\underset{\left( 00\right) }{\overset{H}{\Delta }}\overset{}{%
_{\beta \gamma }^{\alpha }},\underset{\left( 10\right) }{\overset{V}{\Delta }%
}\overset{}{_{\beta \gamma }^{\alpha }}$ \textit{are given by} \textit{(\ref%
{deltaHV}).}\newline
\textbf{Proposition 3.5} \textit{The torsion tensors of the N-linear
connections }\r{D}$\Gamma \left( \mathring{N}\right) $\textit{, }D$^{\top
}\Gamma \left( \check{N}\right) $\textit{\ are related by:}%
\begin{equation*}
\begin{array}{rcl}
\mathring{T}\left( \check{X}^{\mathring{H}},\check{Y}^{\mathring{H}}\right)
& = & \overset{\top }{T}\left( \check{X}^{\mathring{H}},\check{Y}^{\mathring{%
H}}\right) +\overset{\top }{D}\left( \check{X}^{\mathring{H}},\check{Y}^{%
\mathring{H}}\right) -\overset{\top }{D}\left( \check{Y}^{\mathring{H}},%
\check{X}^{\mathring{H}}\right) \\ 
&  &  \\ 
\mathring{T}\left( \check{X}^{\mathring{H}},\check{Y}^{\mathring{V}}\right)
& = & \overset{\top }{T}\left( \check{X}^{\mathring{H}},\check{Y}^{\mathring{%
V}}\right) +\overset{\top }{D}\left( \check{X}^{\mathring{H}},\check{Y}^{%
\mathring{V}}\right) \\ 
&  &  \\ 
\mathring{T}\left( \check{X}^{\mathring{V}},\check{Y}^{\mathring{V}}\right)
& = & \overset{\top }{T}\left( \check{X}^{\mathring{V}},\check{Y}^{\mathring{%
V}}\right), \forall \check{X},\check{Y}\in \mathcal{\check{X}}\left( \widetilde{Osc%
\check{M}}\right).%
\end{array}%
\end{equation*}%
\newline
\textbf{Proposition 3.6} \textit{The torsion d-tensors of the induced
tangent connection of the Cartan metrical }$N$\textit{-linear connection }$%
D\Gamma \left( N\right) $\textit{\ and of the intrinsic Cartan metrical }$N$%
\textit{-linear connection }$\mathring{D}\Gamma \left( \mathring{N}\right) $%
\textit{\ are related by}: 
\begin{equation*}
\begin{array}{lll}
\underset{\left( 00\right) }{\overset{H}{\mathring{T}}}\overset{}{_{\beta
\gamma }^{\alpha }}=\underset{\left( 00\right) }{\overset{H}{T}}\overset{}{%
_{\beta \gamma }^{\alpha }=0}, &  & \underset{\left( 01\right) }{\overset{V}{%
\mathring{T}}}\overset{}{_{\beta \gamma }^{\alpha }}=\underset{\left(
01\right) }{\overset{V}{T}}\overset{}{_{\beta \gamma }^{\alpha }}+\overset{}{%
\underset{00}{\overset{1}{D}}}\underset{}{_{\beta \gamma }^{\alpha }}, \\ 
&  &  \\ 
\underset{\left( 10\right) }{\overset{H}{\mathring{P}}}\overset{}{_{\beta
\gamma }^{\alpha }}=\underset{\left( 10\right) }{\overset{H}{P}}\overset{}{%
_{\beta \gamma }^{\alpha }}, &  & \underset{\left( 11\right) }{\overset{V}{%
\mathring{P}}}\overset{}{_{\beta \gamma }^{\alpha }}=\underset{\left(
11\right) }{\overset{V}{P}}\overset{}{_{\beta \gamma }^{\alpha }}+\overset{}{%
\underset{01}{\overset{1}{D}}}\underset{}{_{\beta \gamma }^{\alpha }}-%
\underset{\left( 10\right) }{\overset{V}{\Delta }}\underset{}{_{\gamma \beta
}^{\alpha }}, \\ 
&  &  \\ 
\underset{\left( 11\right) }{\overset{V}{\mathring{S}}}\overset{}{_{\beta
\gamma }^{\alpha }}=\underset{\left( 11\right) }{S}\overset{}{_{\beta \gamma
}^{\alpha }}=0, &  & 
\end{array}%
\end{equation*}%
\textit{where} $\overset{}{\underset{00}{\overset{1}{D}}}\underset{}{_{\beta
\gamma }^{\alpha }},\overset{}{\underset{01}{\overset{1}{D}}}\underset{}{%
_{\beta \gamma }^{\alpha }},\underset{\left( 10\right) }{\overset{V}{\Delta }%
}\overset{}{_{\beta \gamma }^{\alpha }}$ \textit{are given by} \textit{(\ref%
{D100}), (\ref{D101})\ and (\ref{deltaHV}).
\newline
}\textbf{Proposition 3.7} \textit{The curvature 2-forms }$\mathring{R}$%
\textit{\ and }$\overset{\top }{R}$\textit{\ of the linear connections \r{D%
}}$\Gamma \left( \mathring{N}\right) $\textit{\ and D}$^{\top }\Gamma \left( 
\check{N}\right) $\textit{\ are related by:}%
\begin{equation*}
\begin{array}{ccl}
\mathring{R}\left( \check{X}^{\mathring{H}},\check{Y}^{\mathring{H}}\right) 
\check{Z}^{\mathring{H}} & = & \overset{\top }{R}\left( \check{X}^{\mathring{%
H}},\check{Y}^{\mathring{H}}\right) \check{Z}^{\mathring{H}}+ \\ 
&  &  \\ 
&  & +\left( \mathring{D}_{\check{X}^{\mathring{H}}}\overset{\top }{D}%
\right) \left( \check{Y}^{\mathring{H}},\check{Z}^{\mathring{H}}\right)
-\left( \mathring{D}_{\check{Y}^{\mathring{H}}}\overset{\top }{D}\right)
\left( \check{X}^{\mathring{H}},\check{Z}^{\mathring{H}}\right) \\ 
&  &  \\ 
&  & +\overset{\top }{D}\left( \check{Y}^{\mathring{H}},\overset{\top }{D}%
\left( \check{X}^{\mathring{H}},\check{Z}^{\mathring{H}}\right) \right) -%
\overset{\top }{D}\left( \check{X}^{\mathring{H}},\overset{\top }{D}\left( 
\check{Y}^{\mathring{H}},\check{Z}^{\mathring{H}}\right) \right) \\ 
&  &  \\ 
&  & +\overset{\top }{D}\left( \mathring{D}_{\check{X}^{\mathring{H}}}\check{%
Y}^{\mathring{H}},\check{Z}^{\mathring{H}}\right) -\overset{\top }{D}\left( 
\mathring{D}_{\check{Y}^{\mathring{H}}}\check{X}^{\mathring{H}},\check{Z}^{%
\mathring{H}}\right) ,%
\end{array}%
\end{equation*}%
\begin{equation*}
\begin{array}{l}
\begin{array}{ccl}
\mathring{R}\left( \check{X}^{\mathring{H}},\check{Y}^{\mathring{H}}\right) 
\check{Z}^{\mathring{V}} & = & \overset{\top }{R}\left( \check{X}^{\mathring{%
H}},\check{Y}^{\mathring{H}}\right) \check{Z}^{\mathring{V}}+ \\ 
&  &  \\ 
&  & +\left( \mathring{D}_{\check{X}^{\mathring{H}}}\overset{\top }{D}%
\right) \left( \check{Y}^{\mathring{H}},\check{Z}^{\mathring{V}}\right)
-\left( \mathring{D}_{\check{Y}^{\mathring{H}}}\overset{\top }{D}\right)
\left( \check{X}^{\mathring{H}},\check{Z}^{\mathring{V}}\right) \\ 
&  &  \\ 
&  & +\overset{\top }{D}\left( \check{Y}^{\mathring{H}},\overset{\top }{D}%
\left( \check{X}^{\mathring{H}},\check{Z}^{\mathring{V}}\right) \right) -%
\overset{\top }{D}\left( \check{X}^{\mathring{H}},\overset{\top }{D}\left( 
\check{Y}^{\mathring{H}},\check{Z}^{\mathring{V}}\right) \right) \\ 
&  &  \\ 
&  & +\overset{\top }{D}\left( \mathring{D}_{\check{X}^{\mathring{H}}}\check{%
Y}^{\mathring{H}},\check{Z}^{\mathring{V}}\right) -\overset{\top }{D}\left( 
\mathring{D}_{\check{Y}^{\mathring{H}}}\check{X}^{\mathring{H}},\check{Z}^{%
\mathring{V}}\right) ,%
\end{array}
\\ 
\\ 
\begin{array}{ccl}
\mathring{R}\left( \check{X}^{\mathring{V}},\check{Y}^{\mathring{H}}\right) 
\check{Z}^{\mathring{H}} & = & \overset{\top }{R}\left( \check{X}^{\mathring{%
V}},\check{Y}^{\mathring{H}}\right) \check{Z}^{\mathring{H}}+ \\ 
&  &  \\ 
&  & +\left( \mathring{D}_{\check{X}^{\mathring{V}}}\overset{\top }{D}%
\right) \left( \check{Y}^{\mathring{H}},\check{Z}^{\mathring{H}}\right) +%
\overset{\top }{D}\left( \mathring{D}_{\check{X}^{\mathring{V}}}\check{Y}^{%
\mathring{H}},\check{Z}^{\mathring{H}}\right) ,%
\end{array}
\\ 
\\ 
\begin{array}{ccl}
\mathring{R}\left( \check{X}^{\mathring{V}},\check{Y}^{\mathring{H}}\right) 
\check{Z}^{\mathring{V}} & = & \overset{\top }{R}\left( \check{X}^{\mathring{%
V}},\check{Y}^{\mathring{H}}\right) \check{Z}^{\mathring{V}}+ \\ 
&  &  \\ 
&  & +\left( \mathring{D}_{\check{X}^{\mathring{V}}}\overset{\top }{D}%
\right) \left( \check{Y}^{\mathring{H}},\check{Z}^{\mathring{V}}\right) +%
\overset{\top }{D}\left( \mathring{D}_{\check{X}^{\mathring{V}}}\check{Y}^{%
\mathring{H}},\check{Z}^{\mathring{V}}\right) ,%
\end{array}
\\ 
\\ 
\begin{array}{ccc}
\mathring{R}\left( \check{X}^{\mathring{V}},\check{Y}^{\mathring{V}}\right) 
\check{Z}^{\mathring{H}} & = & \overset{\top }{R}\left( \check{X}^{\mathring{%
V}},\check{Y}^{\mathring{V}}\right) \check{Z}^{\mathring{H}},%
\end{array}
\\ 
\\ 
\begin{array}{ccc}
\mathring{R}\left( \check{X}^{\mathring{V}},\check{Y}^{\mathring{V}}\right) 
\check{Z}^{\mathring{V}} & = & \overset{\top }{R}\left( \check{X}^{\mathring{%
V}},\check{Y}^{\mathring{V}}\right) \check{Z}^{\mathring{V}}, \forall \check{X},\check{Y},\check{Z}\in \mathcal{\check{X}}\left( 
\widetilde{Osc\check{M}}\right) .%
\end{array}%
\end{array}%
\end{equation*}%
\newline
\textbf{Proposition 3.8} \textit{The curvature d-tensors of the induced
tangent connection of the Cartan metrical $N$-linear connection $D\Gamma
\left( N\right) $ and of the intrinsic Cartan metrical $N$-linear connection 
}$\mathring{D}\Gamma \left( \mathring{N}\right) $\textit{\ of the
submanifold }$\widetilde{Osc\check{M}}$\textit{\ are related by}:%
\begin{equation*}
\begin{array}{lll}
\underset{\left( 00\right) }{\overset{H}{\mathring{R}}}\underset{}{_{\beta }}%
\overset{}{^{\alpha }}_{\gamma \delta }=\underset{\left( 00\right) }{\overset%
{H}{R}}\underset{}{_{\beta }}\overset{}{^{\alpha }}_{\gamma \delta }+%
\underset{\left( 00\right) }{\overset{H}{\Delta }}\underset{}{_{\beta }}%
\overset{}{^{\alpha }}_{\gamma \delta }, &  & \underset{\left( 01\right) }{%
\overset{V}{\mathring{R}}}\underset{}{_{\beta }}\overset{}{^{\alpha }}%
_{\gamma \delta }=\underset{\left( 01\right) }{\overset{V}{R}}\underset{}{%
_{\beta }}\overset{}{^{\alpha }}_{\gamma \delta }+\underset{\left( 01\right) 
}{\overset{V}{\Delta }}\underset{}{_{\beta }}\overset{}{^{\alpha }}_{\gamma
\delta }, \\ 
&  &  \\ 
\underset{\left( 10\right) }{\overset{H}{\mathring{P}}}\underset{}{_{\beta }}%
\overset{}{^{\alpha }}_{\gamma \delta }=\underset{\left( 10\right) }{\overset%
{H}{P}}\underset{}{_{\beta }}\overset{}{^{\alpha }}_{\gamma \delta }+%
\underset{\left( 10\right) }{\overset{H}{\Delta }}\underset{}{_{\beta }}%
\overset{}{^{\alpha }}_{\gamma \delta }, &  & \underset{\left( 11\right) }{%
\overset{V}{\mathring{P}}}\underset{}{_{\beta }}\overset{}{^{\alpha }}%
_{\gamma \delta }=\underset{\left( 11\right) }{\overset{V}{P}}\underset{}{%
_{\beta }}\overset{}{^{\alpha }}_{\gamma \delta }+\underset{\left( 11\right) 
}{\overset{V}{\Delta }}\underset{}{_{\beta }}\overset{}{^{\alpha }}_{\gamma
\delta }, \\ 
&  &  \\ 
\underset{\left( 1i\right) }{\overset{V_{i}}{\mathring{S}}}\underset{}{%
_{\beta }}\overset{}{^{\alpha }}_{\gamma \delta }=\underset{\left( 1i\right) 
}{\overset{V_{i}}{S}}\underset{}{_{\beta }}\overset{}{^{\alpha }}_{\gamma
\delta }\text{, }\left( i=0,1;V_{0}=H\right) , &  & 
\end{array}%
\end{equation*}%
\newline
\textit{where}

\begin{equation}
\begin{array}{lll}
\underset{\left( 00\right) }{\overset{H}{\Delta }}\underset{}{_{\beta }}%
\overset{}{^{\alpha }}_{\gamma \delta } & = & \delta _{\delta }\underset{%
\left( 00\right) }{\overset{H}{\Delta }}\underset{}{_{\beta \gamma }^{\alpha
}}-\delta _{\gamma }\underset{\left( 00\right) }{\overset{H}{\Delta }}%
\underset{}{_{\beta \delta }^{\alpha }}-D^{\varepsilon }\overset{}{_{\delta }%
}\dot{\partial}_{1\varepsilon }\underset{\left( 00\right) }{\overset{H}{L}}%
\overset{}{_{\beta \gamma }^{\alpha }}-D^{\varepsilon }\overset{}{_{\delta }}%
\dot{\partial}_{1\varepsilon }\underset{\left( 00\right) }{\overset{H}{%
\Delta }}\underset{}{_{\beta \gamma }^{\alpha }} \\ 
&  &  \\ 
&  & +D^{\varepsilon }\overset{}{_{\gamma }}\dot{\partial}_{1\varepsilon }%
\underset{\left( 00\right) }{\overset{H}{L}}\overset{}{_{\beta \delta
}^{\alpha }}+D^{\varepsilon }\overset{}{_{\gamma }}\dot{\partial}%
_{1\varepsilon }\underset{\left( 00\right) }{\overset{H}{\Delta }}\underset{}%
{_{\beta \delta }^{\alpha }} \\ 
&  &  \\ 
&  & +\underset{\left( 00\right) }{\overset{H}{L}}\overset{}{_{\beta \gamma
}^{\sigma }}\underset{\left( 00\right) }{\overset{H}{\Delta }}\underset{}{%
_{\sigma \delta }^{\alpha }}+\underset{\left( 00\right) }{\overset{H}{\Delta 
}}\underset{}{_{\beta \gamma }^{\sigma }}\underset{\left( 00\right) }{%
\overset{H}{L}}\overset{}{_{\sigma \delta }^{\alpha }}-\underset{\left(
00\right) }{\overset{H}{L}}\overset{}{_{\beta \delta }^{\sigma }}\underset{%
\left( 00\right) }{\overset{H}{\Delta }}\underset{}{_{\sigma \gamma
}^{\alpha }}-\underset{\left( 00\right) }{\overset{H}{\Delta }}\underset{}{%
_{\beta \delta }^{\sigma }}\underset{\left( 00\right) }{\overset{H}{L}}%
\overset{}{_{\sigma \gamma }^{\alpha }} \\ 
&  &  \\ 
&  & +\underset{\left( 01\right) }{\overset{H}{C}}\overset{}{_{\beta \sigma
}^{\alpha }}\underset{00}{\overset{1}{D}}\underset{}{_{\beta \gamma
}^{\sigma }},%
\end{array}
\label{deltaH00}
\end{equation}%
\begin{equation}
\begin{array}{lll}
\underset{\left( 01\right) }{\overset{V}{\Delta }}\underset{}{_{\beta }}%
\overset{}{^{\alpha }}_{\gamma \delta } & = & \overset{}{\underset{1}{D}%
^{\varepsilon }}\overset{}{_{\gamma }}\dot{\partial}_{1\varepsilon }\underset%
{\left( 10\right) }{\overset{V}{L}}\overset{}{_{\beta \delta }^{\alpha }}-%
\overset{}{\underset{1}{D}^{\varepsilon }}\overset{}{_{\delta }}\dot{\partial%
}_{1\varepsilon }\underset{\left( 10\right) }{\overset{V}{L}}\overset{}{%
_{\beta \gamma }^{\alpha }}+\dfrac{1}{2}\delta _{\delta }\left( \underset{%
\left( 10\right) }{\overset{V}{\Delta }}\overset{}{_{\beta \gamma }^{\alpha }%
}\right) -\dfrac{1}{2}\delta _{\gamma }\left( \underset{\left( 10\right) }{%
\overset{V}{\Delta }}\overset{}{_{\beta \delta }^{\alpha }}\right) \\ 
&  &  \\ 
&  & +\dfrac{1}{2}\overset{}{\underset{1}{D}^{\varepsilon }}\overset{}{%
_{\gamma }}\dot{\partial}_{1\varepsilon }\left( \underset{\left( 10\right) }{%
\overset{V}{\Delta }}\overset{}{_{\beta \delta }^{\alpha }}\right) -\dfrac{1%
}{2}\overset{}{\underset{1}{D}^{\varepsilon }}\overset{}{_{\delta }}\dot{%
\partial}_{1\varepsilon }\left( \underset{\left( 10\right) }{\overset{V}{%
\Delta }}\overset{}{_{\beta \gamma }^{\alpha }}\right) + \\ 
&  &  \\ 
&  & \dfrac{1}{2}\underset{\left( 10\right) }{\overset{V}{L}}\overset{}{%
_{\beta \gamma }^{\varepsilon }}\underset{\left( 10\right) }{\overset{V}{%
\Delta }}\overset{}{_{\varepsilon \delta }^{\alpha }}-\dfrac{1}{2}\underset{%
\left( 10\right) }{\overset{V}{L}}\overset{}{_{\beta \delta }^{\varepsilon }}%
\underset{\left( 10\right) }{\overset{V}{\Delta }}\overset{}{_{\varepsilon
\gamma }^{\alpha }}+\dfrac{1}{2}\underset{\left( 10\right) }{\overset{V}{%
\Delta }}\overset{}{_{\beta \gamma }^{\varepsilon }}\underset{\left(
10\right) }{\overset{V}{L}}\overset{}{_{\varepsilon \delta }^{\alpha }}-%
\dfrac{1}{2}\underset{\left( 10\right) }{\overset{V}{\Delta }}\overset{}{%
_{\beta \delta }^{\varepsilon }}\underset{\left( 10\right) }{\overset{V}{L}}%
\overset{}{_{\varepsilon \gamma }^{\alpha }} \\ 
&  &  \\ 
&  & +\dfrac{1}{4}\left( \underset{\left( 10\right) }{\overset{V}{\Delta }}%
\overset{}{_{\varepsilon \delta }^{\alpha }}\underset{\left( 10\right) }{%
\overset{V}{\Delta }}\overset{}{_{\beta \gamma }^{\varepsilon }}-\underset{%
\left( 10\right) }{\overset{V}{\Delta }}\overset{}{_{\varepsilon \gamma
}^{\alpha }}\underset{\left( 10\right) }{\overset{V}{\Delta }}\overset{}{%
_{\beta \delta }^{\varepsilon }}\right) +\underset{\left( 11\right) }{%
\overset{V}{C}}\overset{}{_{\beta \sigma }^{\alpha }}\underset{\left(
10\right) }{\overset{V}{\Delta }}\overset{}{_{\gamma \delta }^{\sigma }},%
\end{array}
\label{deltaV01}
\end{equation}%
\begin{equation}
\begin{array}{lll}
\underset{\left( 10\right) }{\overset{V}{\Delta }}\underset{}{_{\beta }}%
\overset{}{^{\alpha }}_{\gamma \delta } & = & -\underset{\left( 00\right) }{%
\overset{H}{\Delta }}\underset{}{_{\varepsilon \gamma }^{\alpha }}\underset{%
\left( 01\right) }{\overset{H}{C}}\overset{}{_{\beta \delta }^{\varepsilon }}%
+\underset{\left( 00\right) }{\overset{H}{\Delta }}\underset{}{_{\beta
\gamma }^{\varepsilon }}\underset{\left( 01\right) }{\overset{H}{C}}\overset{%
}{_{\varepsilon \delta }^{\alpha }}+\underset{\left( 01\right) }{\overset{H}{%
C}}\overset{}{_{\beta \varepsilon }^{\alpha }}\overset{}{\underset{01}{%
\overset{1}{D}}}\underset{}{_{\gamma \delta }^{\varepsilon }}+\dot{\partial}%
_{1\delta }\underset{\left( 00\right) }{\overset{H}{\Delta }}\underset{}{%
_{\beta \gamma }^{\alpha }} \\ 
&  &  \\ 
&  & +D^{\sigma }\overset{}{_{\gamma }}\dot{\partial}_{1\sigma }\underset{%
\left( 01\right) }{\overset{H}{C}}\overset{}{_{\beta \delta }^{\alpha }},%
\end{array}
\label{deltaV10}
\end{equation}%
\begin{equation}
\begin{array}{lll}
\underset{\left( 11\right) }{\overset{V}{\Delta }}\underset{}{_{\beta }}%
\overset{}{^{\alpha }}_{\gamma \delta } & = & \dfrac{1}{2}\dot{\partial}%
_{1\delta }\left( \underset{\left( 10\right) }{\overset{V}{\Delta }}\overset{%
}{_{\beta \gamma }^{\alpha }}\right) -\dfrac{1}{2}\left( \underset{\left(
10\right) }{\overset{V}{\Delta }}\overset{}{_{\varepsilon \gamma }^{\alpha }}%
\underset{\left( 11\right) }{\overset{V}{C}}\overset{}{_{\beta \delta
}^{\varepsilon }}-\underset{\left( 11\right) }{\overset{V}{C}}\overset{}{%
_{\varepsilon \delta }^{\alpha }}\underset{\left( 10\right) }{\overset{V}{%
\Delta }}\overset{}{_{\beta \gamma }^{\varepsilon }}\right. - \\ 
&  &  \\ 
&  & \left. -\underset{\left( 11\right) }{\overset{V}{C}}\overset{}{_{\beta
\varepsilon }^{\alpha }}\left( \underset{\left( 10\right) }{\overset{V}{%
\Delta }}\overset{}{_{\delta \gamma }^{\varepsilon }}-\underset{\left(
10\right) }{\overset{V}{\Delta }}\overset{}{_{\gamma \delta }^{\varepsilon }}%
\right) \right) ,%
\end{array}
\label{deltaV11}
\end{equation}%
\textit{where} $\underset{\left( 00\right) }{\overset{H}{\Delta }}\overset{}{%
_{\beta \gamma }^{\alpha }},\underset{\left( 10\right) }{\overset{V}{\Delta }%
}\overset{}{_{\beta \gamma }^{\alpha }}$ \textit{are given by (\ref{deltaHV}%
)}.\newline

\end{document}